\documentclass{scrartcl}
              \usepackage[utf8]{inputenc}
\usepackage[T1]{fontenc}
\usepackage{lmodern}
\usepackage{fixltx2e}
\usepackage{booktabs}
\usepackage{microtype}
\usepackage{mathtools}
\usepackage{graphicx}
\usepackage{amsfonts}
\usepackage{color}
\usepackage{flafter}
\usepackage{amsmath,amssymb,amsthm,cancel,units}
\usepackage{pgfplots}
\usepackage{algorithmic} \usepackage[ruled]{algorithm}
\usepackage[unicode=true,plainpages=false]{hyperref}
\hypersetup{colorlinks=true,linkcolor=magenta,anchorcolor=magenta,urlcolor=blue,citecolor=blue}

\usetikzlibrary{calc,patterns,decorations.pathmorphing,decorations.markings}
\newtheorem{theorem}{Theorem}[section]
\newtheorem{definition}[theorem]{Definition}

\usepackage[a4paper]{geometry}
\usepackage{graphicx}
\usepackage{longtable}
\usepackage{float}

\title{Fast low-rank solution of the Poisson equation with application to the Stokes problem}
\author{E.A. Muravleva, I.V. Oseledets}
\date{\today}
\hypersetup{
  pdfkeywords={},
  pdfsubject={},
  pdfcreator={Emacs Org-mode version 7.8.11}}

\begin{document}

\maketitle

\section{Introduction}
\label{sec-1}

Solvers for the Poisson equation are important in many different
application areas, thus their improvement is of great interest even
for special grids and right-hand sides. In this paper we consider a
very specific class of Poisson equations: Poisson equations in two and
three dimensions on
tensor-product uniform grids with \emph{low-rank} right-hand sides. In this
case, it is possible to reduce the complexity of the solver to an
almost linear complexity in a one-dimensional grid size.
 The problems of such kind have been considered previously by several authors \cite{grasedyck-kron-2004,
GHK-ten_inverse_ellipt-2005, beylkin-2002} and the linear complexity
is indeed possible. However, in small dimensions (especially in two
dimensions)  the constant hidden in $\mathcal{O}(n)$ can be very high,
and the full representation is more efficient for a wide range of
$n$. The main goal of
this paper is, staying within the framework of low-rank
approximations, provide a new and efficient algorithm for the approximate
solution  the Poisson equation. 

To show the effectiveness of our approach, we use it to create a fast
solver for Stokes problem. The Stokes problem is one of the classical problems of mathematical
physics. It is often encountered as a subproblem in more
complicated problems, such as Navier-Stokes problem, unsteady Stokes
problem, flows of non-Newtonian fluids. Numerical methods for the
solution of the Stokes problem are a classical topic of the
computational fluid dynamics, thus any improvement of the efficiency
of the Stokes solvers (at least in some important cases) is of great
interest.  The Uzawa method for the Stokes problem requires the
repetitive solution of the Poisson equations. The whole iterative procedure
(including matrix-by-vector products, arithmetics, dot products) will
be implemented in the low-rank format, and it is not always an easy
task, since at each step the accuracy should be monitored to avoid the
growth of the number of parameters and the growth of the error.

The method proposed in this paper has its limitations: it works for a
special  discretization of the Laplace operator on
tensor-product grids and for special right-hand sides, but the class
of right-hand sides in two and three dimensions is not small: it
includes sufficiently smooth functions (i.e., approximated by
polynomials), so-called asymptotically smooth functions, sparse
right-hand sides and so on. The approximation work is done on the
algebraic level using unified algorithms. To get better complexity and
accuracy, one can use advanced discretization techniques: $hp$-methods, high-order schemes,
 spectral and pseudospectral approaches, etc. However, they require
 significant additional work on the analytic level.  On the other
 hand, one can use certain approximations on the discrete level, by
 approximating the discrete solution using some low-parametric
 format. As such a format, we will use a low-rank factorization of
 matrices (in two dimensions) and tensors (in three dimensions). The
 approach can be used in arbitrary dimension, but we will leave this
 topic for future work.
\section{Model problem: discretization and solution method}
\label{sec-2}

The Stokes problem is used a basic problem to illustrate the
low-rank techniques that are used.  The numerical scheme consists of
three basic steps.
\begin{enumerate}
\item Take a simple discretization of the Stokes problem (we will use
   semi-staggered grids)
\item Use a mesh-size-independent convergent iterative scheme (we will
   use the Uzawa method)
\item Replace each step of the method by operations in the low-parametric
   format with truncation.
\end{enumerate}

The whole procedure is very simple, and as it will be seen later on,
most of the steps are also simple in the tensor format. However, to get a more efficient
method, several modifications should be made to the ``naive''
approach. For the Stokes problem, the second step requires the
solution of the Poisson problem at each iteration step. Our main
result is a new algorithm for the solution of such problem, based on
the cross approximation \cite{tee-cross-2000, bebe-2000} in the frequency space. The well-known
approaches for the solution of the Poisson equation in low-rank
formats rely on a special approximation to the inverse of the Laplace
operator by a sum of tensor (Kronecker) products, but in two and three
dimensions the complexity can be reduced signicantly using a new approach.
 
\begin{center}
\begin{figure}[!htbp]
\resizebox{10cm}{!}
{
\input{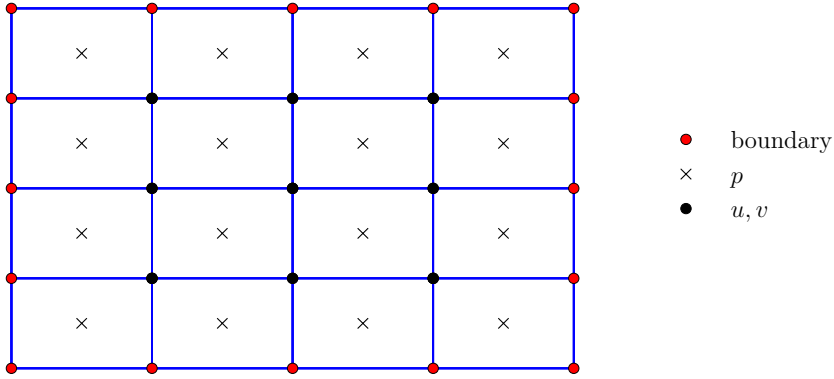}
}
\caption{The model problem: semi-staggered grids for velocity components and pressure}\label{flr:grid}
\end{figure}
\end{center}

We consider the Stokes problem in a rectangular domain $\Omega = [0,1]^d$, $d=2,3$.
  The problem is discretized using the finite difference method (FDM)
  on semi-staggered grids. For details, see \cite{om-stockes-2010,mur-phd-2010,mur-ker-2008}. We will give
  only the final matrix form of the discrete problem. There are $d+1$
  unknown vectors ($d$ for velocity components and one for
  pressure). The components of velocity are defined in the vertices of
  the grid cells, and the pressure is defined in the centers of the
  cells (see Figure \ref{flr:grid}). The grid for pressure is $n_x \times n_y$ in
  two-dimensional case, and $n_x \times n_y \times n_z$ in
  three-dimensional case, respectively. For simplicity we will assume
  that the grid sizes are equal: $n_x = n_y = n_z = n$. The mesh size
  is defined as $h = \frac{1}{n}$.

The discretized Stokes problem has the following form:
 \begin{equation}\label{flr:stprob}
 \begin{pmatrix}
A & B\\
B^T & 0
\end{pmatrix}
\begin{pmatrix}
u \\ p
\end{pmatrix} = 
\begin{pmatrix}
f \\ g
\end{pmatrix}.
\end{equation}
We use the standard notations: $u$ is the velocity vector, $p$ is the
pressure.
Here $A$ is a $d \times d$ block diagonal matrix of the form

$$ A = I_d \otimes \Delta, $$
where $\Delta$ is a discrete Laplace operator which will be defined
later on, and $I_d$ is an identity matrix of order $d$.

The matrix $B$ is a block matrix matrix of the form

$$
    B = \begin{pmatrix}
B_{x}^{\top} &  B_{y}^{\top}
\end{pmatrix}^{\top}
$$
in two-dimensional case, and 
$$
       B = \begin{pmatrix}
B_{x}^{\top} &  B_{y}^{\top} & B_{z}^{\top}
\end{pmatrix}^{\top}
$$
in three-dimensional case. 

To define the components of the discrete gradient operators, it is
convenient to use the following auxiliary matrices, $G$ and $H$.

They are defined as

$$ G = E - Z, \quad H  = E + Z, $$
where $E,Z$ are $(n-1) \times n$ matrices. The matrix $E$ deletes the
first element of the vector (and moves all others up by one), and the matrix $Z$ deletes the last
element of the vector (and moves all others down by one).

In two-dimensional case, the matrices $B_x$ and $B_y$ have the form

$$
   B_x = \frac{1}{4h} G \otimes H, \quad B_y = \frac{1}{4h} H \otimes G.
$$
In three-dimensional case, the matrices $B_x$, $B_y$ and $B_z$ have
the form:

$$
   B_x = \frac{1}{4h} G \otimes H \otimes H, \quad B_y = \frac{1}{4h} H
   \otimes G \otimes H, \quad B_z = \frac{1}{4h} H \otimes H \otimes G.
$$

Finally, the discrete Laplace operator has the following form:
\begin{equation*}
\begin{split}
 &  \Delta = B_x B^{\top}_x + B_y B^{\top}_y, \\
 & \Delta = B_x B^{\top}_x + B_y B^{\top}_y + B_z B^{\top}_z, 
\end{split}
\end{equation*}
in two and three dimensions respectively.

Let us note that this is a non-standard 9-point stencil in 2D
(27-point stencil in 3D) discretization for the Laplace operator, but
its main benefit is that it is \emph{consistent} with the discrete gradient operator.

  The standard approach to solve \eqref{flr:stprob} is the Uzawa
  method \cite{benzi-saddle-2005}. The system is reduces to the following equation
  for the pressure:

\begin{equation}\label{flr:schur}
(B^{\top} A^{-1} B) p = B^{\top} A^{-1} f.
\end{equation}

The matrix of the system \eqref{flr:schur} is positive semi-definite,
and the conjugate-gradient (CG) method is a convenient tool to solve
it. An important question is the spectrum of the Schur operator. It is
known \cite{om-stockes-2010, mur-phd-2010} that it has $2$ zero eigenvalues in 2D, and
$3n-1$ eigenvalues in $3d$, and $(n-2)^d$ eigenvalues, equal to $1$. 
The numerical experiments confirm, that the remaining part of the
spectrum lies on $[0,1]$, and is bounded from below by a constant,
independent of $h$. The linear system
\eqref{flr:schur} is consistent, i.e., the right-hand side is
orthogonal to any vector $q$ in the kernel of
the Schur complement, thus the conjugate gradient method effectively
works on the subspace, orthogonal to the kernel, where the operator is
well-conditioned, thus the total number of iterations does not depend on
$h$. 
\section{Low-rank formats}
\label{sec-3}

 Our main goal is the solution of the problem \eqref{flr:schur} in
 low-rank formats. In two dimensions, the discrete pressure is a
 vector of length $(n-1)^2$, and it can be naturally considered as an
 $(n-1) \times (n-1)$ matrix. In a three-dimensional case, the discrete
 pressure can be naturally considered as an $(n-1) \times (n-1) \times
 (n-1)$ three-dimensional tensor. We want to approximate those
 tensors, effectively reducing the number of parameters and the
 computational cost, while maintaining the required accuracy
 $\varepsilon$ of the computations. This accuracy should be consistent
 with the mesh discretization error, which is $\mathcal{O}(h^2)$ for
 the pressure.

Let us recall some basic facts about low-rank tensor
 decompositions in two and three dimensions. 
\begin{definition}
  The $n_1 \times n_2$ matrix $M$ is said to be in the low-rank format
  with rank $r$, if it can be represented as
 \begin{equation}\label{flr:lr}
    M = U V^{\top} = \sum_{\alpha=1}^r u_{\alpha} v_{\alpha}^{\top},
 \end{equation}
where $U$ is an $n_1 \times r$ matrix, and $V$ is an $n_2 \times r$
matrix.
\end{definition}
 The dyadic decomposition \eqref{flr:lr} can be computed via
the singular value decomposition (SVD). The linear operator, acting on
 a space of $n_1 \times n_2$ matrices, can be represented as a $(n_1
 n_2) \times (n_1 n_2)$ matrix. The corresponding low-rank format for
 the matrix is given in the following form.
\begin{definition}
  An $(n_1 n_2) \times (n_1 n_2)$ matrix $A$ is said to be in the
  low-rank (Kronecker) format with rank $R$, if it can be represented
  as a sum of Kronecker products:
\begin{equation}\label{flr:lrmat}
  A = \sum_{\alpha=1}^R C_{\alpha} \otimes D_{\alpha},
\end{equation}
where $C_{\alpha}$ is an $n_1 \times n_1$ matrix, and $D_{\alpha}$ is
an $n_2 \times n_2$ matrix, and $\otimes$ is a Kronecker product of matrices.
\end{definition}

The low-rank format for operators and vectors allows for fast linear
algebra operations
\cite{ost-latensor-2009,ot-hyper-2005,hkt-iter-2008}. For example, 
multiplication of a matrix of tensor rank $R$ by a vector of rank $r$ yields a vector of rank $Rr$.

An important operation is the so-called \emph{rounding}. Most basic
operations (like addition or matrix-by-vector product) increase the
rank, thus the reapproximation with some specified accuracy is
required. If the matrix is given in the low-rank format, then the
rounding procedure can be implemented without the computation of the
full SVD. Its complexity with respect to $n,m,r$ is $\mathcal{O}(
(n+m) r^2 + r^3)$, i.e., it is linear in the mode size (compared to
cubic for the full SVD). 

In three dimensions, the so-called Tucker format can be very useful:
\cite{Tucker,lathauwer-svd-2000,khor-ml-2009}. 

\begin{definition}
An $n_1 \times n_2 \times n_3$ three-dimensional tensor $A$ is said to be
in the Tucker format, if it can be represented as
\begin{equation}\label{flr:tucker}
  A(i_1,i_2,i_3) = \sum_{\alpha_1,\alpha_2,\alpha_3}
  G(\alpha_1,\alpha_2,\alpha_3) U_1(i_1,\alpha_1) U_2(i_2,\alpha_2) U_3(i_3,\alpha_3),
\end{equation}
where the numbers $\alpha_k$ vary from $1$ to $r_k$.

The matrices $U_k = [U_k(i_k,\alpha_k)]$ are called \emph{factors} of the
Tucker decomposition, the numbers $r_k$ are called \emph{Tucker
ranks}, and the three-dimensional tensor $G$ is called the
\emph{Tucker core}.
\end{definition}

The Tucker format for the matrix is defined in the same fashion. The
basic linear operations for the Tucker format can also be defined, as
well as a fast rounding procedure with complexity that is linear in
the mode size \cite{ost-latensor-2009}.
\section{Fast inversion of the Laplacian in tensor formats}
\label{sec-4}
\subsection{Known approaches}
\label{sec-4-1}

The basic computational core of the algorithm is the solution of the
Poisson equation (see \eqref{flr:schur}): 
\begin{equation}\label{flr:pois}
   \Delta u = g.
\end{equation}
We will consider first the two-dimensional case. The generalization to
higher dimensions will be outlined further. The matrix $\Delta$ in our
case has the form  (up to a scaling factor, which is not important)
$$
    \Delta = A_1 \otimes A_2 + A_2 \otimes A_1,
$$
where
$$
  A_1 = GG^{\top}, A_2 = HH^{\top}.
$$
It is not difficult to see, that the matrices $A_1$ and $A_2$ commute
and are diagonalized by the discrete sine transform (DST):
$$
   A_1 = S \Lambda_1 S^{\top}, A_2 = S \Lambda_2 S^{\top},  
$$ 
thus the full matrix $\Delta$ can be written as
\begin{equation}\label{flr:diag}
    \Delta = (S \otimes S) (\Lambda_1 \otimes \Lambda_2 + \Lambda_2
    \otimes \Lambda_1) (S^{\top} \otimes S^{\top}).
\end{equation}
The representation \eqref{flr:diag} is a textbook way to solve such
kind of problems in the full format, by applying the DST, inverting
the diagonal matrix, and then applying the DST back. 

For the low-rank formats, the situation is as follows. The application
of the DST can be done in a fast way by computing the DST of the factors
 $U$ and $V$, and it is an $\mathcal{O}(nr \log n)$
complexity. Moreover, then the rank does not change. But the inversion of the diagonal matrix
\begin{equation}\label{flr:diagmat}
   \Lambda = \Lambda_1 \otimes \Lambda_2 + \Lambda_2 \otimes \Lambda_1,
\end{equation}
suddenly becomes the bottleneck, since it is a $\mathcal{O}(n^2)$ computation!
Thus, the matrix $\Lambda^{-1} \widehat{g}$ should be approximated in the
low-rank format.

A well-established approach to compute the inverse of $\Lambda^{-1}$    
is based on the approximation by exponential sums
\cite{rokhlin-gauss-1998,GHK-ten_inverse_ellipt-2005,grasedyck-kron-2004,khor-rstruct-2006,hackbra-expsum-2005}. 
It is uses the following quadrature formula:
\begin{equation}\label{flr:int}
   \frac{1}{x} = \int_{0}^{\infty} e^{-px} dp \approx
   \sum_{\alpha=1}^R w_k e^{-p_k x}, 
\end{equation}
where $p_k$ and $w_k$ are quadrature weights. Typical value of $r$
required to achieve a good accuracy is of order of several tens.

The equation \eqref{flr:int} can be used to get the approximation of the
inverse matrix $\Lambda$:
\begin{equation}\label{flr:expsum}
\Lambda^{-1} \approx \sum_{\alpha=1}^R w_k e^{-p_k \Lambda} =
\sum_{\alpha=1}^R w_k e^{-p_k(\Lambda_1 + \Lambda_2)} \otimes e^{-p_k(\Lambda_1 + \Lambda_2)}.
\end{equation}
\subsection{A new method based on cross approximation}
\label{sec-4-2}
  Let us estimate the typical complexity of the method, discussed in
  the previous subsection. The main computational task is the
  multiplication by $\Lambda^{-1}$. Let $r$ be the rank of the vector, and $R$
  be the tensor rank of the inverse of $\Lambda$. Then the result
  would have the rank $Rr$. Thus, in two dimensions, the method will be effective only
  when $n \geq (Rr) \sim 1000$, since the typical values of $R$ and
  $r$ are around $30$. 

  The main cost is due to the approximation of the full inverse
  matrix. However, we need only to multiply this inverse by a vector,
  and the result is expected to have a small rank. This structure is
  not used in the methods based on the approximation of the
  inverse matrix.

  To obtain a new method, let us look more closely on the main
  operation:
  $$ 
      \widehat{f} = \Lambda^{-1} \widehat{g},
  $$ 
  or elementwise:
  $$
      \widehat{f}_{ij} = \frac{\widehat{g}_{ij}}{\Lambda_{ij}}.
  $$
 In our case, 
$$\Lambda_{ij} = (4 - \lambda_i) \lambda_i + (4 -
  \lambda_j) \lambda_j = \mu_i + \mu_j,$$
therefore
\begin{equation}\label{flr:div}
   \widehat{f}_{ij} = \frac{\widetilde{g}_{ij}}{\mu_i + \mu_j}.
\end{equation}

The most important observation is that \eqref{flr:div} is an
\emph{elementwise division} of a matrix of rank $r$ by a matrix of
rank $2$, and we can compute any element of the matrix $\widehat{F} =
[f_{ij}]$ in $\mathcal{O}(r)$ operations. 

Our main assumption is that the result of the division can be in turn
approximated by a matrix of a small rank, $\mathcal{O}(r)$. Thus, we can try to
approximate the result directly, without the approximation of the
$\Lambda^{-1}$. Since we can compute any prescribed element fastly, it
is natural to use the cross method for the approximation of low-rank
matrices (in two dimensions) \cite{tee-cross-2000, bebe-2000} and for
tensors \cite{ost-tucker-2008,ot-ttcross-2010}. In the matrix case, to
approximate a matrix with rank $r$ it is sufficient to compute $r$
rows and $r$ columns of it, using a certain heuristical technique
(a quasioptimal choice can be based on the maximal-volume submatrix
\cite{gt-maxvol-2001} using the maxvol-algorithm
\cite{gostz-maxvol-2010}). Thus, the expected complexity is
$\mathcal{O}(nr^2)$ operations, which is much less, than the
approximation of the full inverse.  
\section{Stokes solver in the low-rank format}
\label{sec-5}

The Uzawa method requires the solution of the  Poisson equation and
multiplication by the gradient matrices, and some addition. All these
operations are done approximately with some accuracy $\varepsilon$. To
solve the linear system with the Schur complement and inexact
matrix-by-vector products, we used the inexact GMRES as implemented in 
\cite{dc-tt_gmres-2013}, which implements the adaptive accuracy for
the intermediate matrix-by-vector products. Despite the fact that our
matrix is positive semidefinite, we decided to use GMRES over the
conjugate gradients in this case, since the GMRES method proved to be
much more robust to the approximate inversion of the Laplacian. 
\section{Numerical experiments}
\label{sec-6}

As a  numerical test, we have selected the following analytical
solution:

\begin{equation}\label{flr:ansin2d}
 \begin{split}
  &u = \frac{1}{4 \pi^2} \sin 2 \pi x (1 - \cos 2 \pi y), \\
  &v = \frac{1}{4 \pi^2} (1 - \cos 2 \pi x) \sin 2 \pi y, \\
  &p = \frac{1}{\pi} \sin 2 \pi x \sin 2 \pi y. 
 \end{split}
\end{equation}
The corresponding right-hand side is then 
\begin{equation}
\begin{split}\label{flr:ansin2d-2}
   &f_u = \sin 2 \pi y, \\
   &f_v = \sin 2 \pi x, \\
   &g = -\frac{1}{\pi}\sin 2 \pi x \sin 2 \pi y.
\end{split}
\end{equation}

The solution for $p$ has a perfect structure (rank-1). However, during
the iterative process the intermediate ranks may grow quite a lot, and
that increases the complexity. The Table \ref{flr:table-1} provides
the dependence of the total computation time for full low-rank method
and for the same method in the full format. The third column gives an
error of the approximation with respect to the true solution. For all
computations the same threshold $\varepsilon = 5 \cdot 10^{-9}$ was
used that is safe to match the discretization error. It is clear, that
the low-rank method is slower for small mode sizes (until $n \leq
256$) and is becoming faster due to better scaling in $n$. The
almost linear scaling in $n$ is also observed.

\begin{table}[H]
\caption{Numerical results for the two-dimensional sine example} \label{flr:table-1}
\begin{center}
\begin{tabular}{rll}
    n  &  Time (LR/full)    &  Rel. error in $p$ (LR/full)  \\
\hline
   64  &  0.17/0.024        &  1.2e-03/1.2e-03              \\
  128  &  0.29/0.061        &  2.9e-04/2.9e-04              \\
  256  &  0.54/0.39         &  7.8e-05/7.8e-05              \\
  512  &  \textbf{1.2/2.3}  &  1.9e-05/1.9e-05              \\
 1024  &  \textbf{2.8/10}   &  8.6e-06/8.2e-06              \\
\end{tabular}
\end{center}
\end{table}

\begin{figure}[!htbp]
\begin{center}
\resizebox{10cm}{!}
{
\input{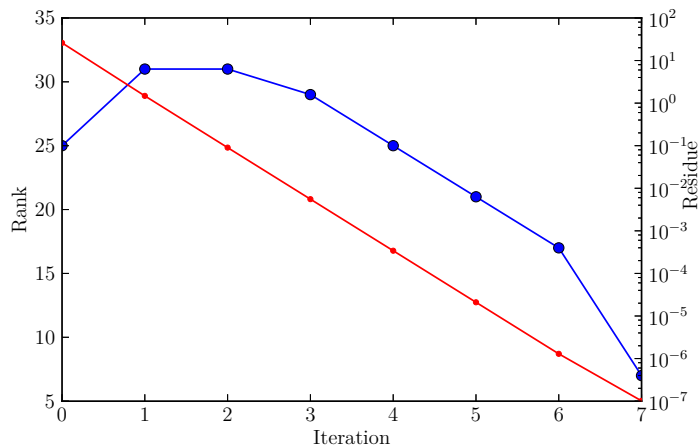}
}
\end{center}
\caption{The ranks of Krylov vectors and the convergence of the
Uzawa-GMRES method}\label{flr:rank-and-res}
\end{figure}

In Figure \ref{flr:rank-and-res} we present the approximate rank of
the Krylov vectors in
the inexact GMRES and the residual computed in the GMRES method. We
see that the rank decreases with the iteration number (which is
consistent with the inexact GMRES theory, the next Krylov vectors are
approximated with less accuracy), and the GMRES converges linearly.

\subsection{Lid-driven cavity flow}
\label{sec-6-1}

\begin{figure}[H]
\begin{tikzpicture}

    \begin{scope}[x=12em, y=12em]
    \tikzstyle{ground}=[fill,pattern=north east lines,draw=none]
     %\draw (-1,0) -- (0,0);
      %\draw (0,0) -- (0,-1);
      %\draw (0,-1) -- (1, -1);
      %\draw (1,-1) -- (1, 0);
      %\draw (1, 0) -- (2, 0);
      \draw[densely dashed, gray] (0,0) -- (1,0);
      \draw[->, gray] (0.25,0.1) -- (0.75, 0.1);
      \node at (0.5, 0.15) {$u = 1, v = 0$};
      \node at (-0.3, -0.5) {$u = 0, v = 0$};
      \node at (1.3, -0.5) {$u = 0, v = 0$};
      \node at (0.5, -1.13) {$u = 0, v = 0$};
      \node (ground1) at (0.5,-1) [anchor = north, ground, minimum width=12em] {};
      \draw (ground1.north west) -- (ground1.north east);
      \node (ground2) at (0, -0.5) [anchor = east, ground, minimum height=12em] {};
      \draw (ground2.south east) -- (ground2.north east);
      \node (ground3) at (1,-0.5) [anchor = west, ground, minimum height=12em] {};
      \draw (ground3.south west) -- (ground3.north west);
      \node (ground4) at (1.5, 0) [anchor = north, ground, minimum width=12em] {};
      \draw (ground4.north west) -- (ground4.north east);
      \node (ground5) at (-0.5, 0) [anchor = north, ground, minimum width=12em] {};
      \draw (ground5.north west) -- (ground5.north east);
  \end{scope}
    % include your tikz code here
\end{tikzpicture}
\caption{The lid-driven cavity flow}\label{flr:lid}
\end{figure}
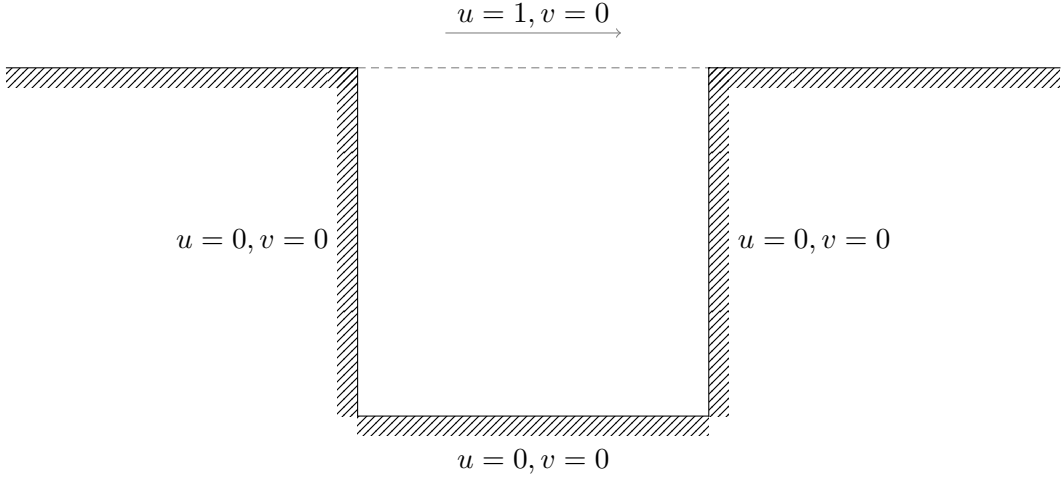

 The second example is the classical lid-driven cavity test (see the
 picture). The homogenius Dirichlet boundary conditions are introduced
 on all the boundaries for $u$ and $v$ except the upper boundary which
 is moving with constant velocity (see Figure \ref{flr:lid}).

The results are presented in Table \ref{flr:table-2}. Since
 now we do not the analytic solution, we can compare the solution
 obtained from the full format with the solution obtained the the
 low-rank method. The sublinear scaling in the number of unknowns is
 clearly visible.
\begin{table}[H]
\caption{Numerical results for the two-dimensional lid-driven cavity flow} \label{flr:table-2}
\begin{center}
\begin{tabular}{rlr}
    n  &  Time (LR/full)    &  Rel. error  \\
  256  &  0.89/0.63         &     6.7e-07  \\
  512  &  \textbf{2.5/3.2}  &     5.6e-07  \\
 1024  &  \textbf{6.3/13}   &     5.9e-07  \\
\end{tabular}
\end{center}
\end{table}

Finally, we will compare the new method with the exponential sums
approach. We have taken $n = 1024$ and $\varepsilon = 1e-7$. That gave
the rank of the inverse $35$. The typical rank of the Krylov vector for the lid-driven
cavity test was around $30$ as well. The timing for the Hadamard
product and then rounding was \textbf{1.6} seconds, and for the cross method
the time was \textbf{0.07} seconds, i.e, more than 20 times faster.
\section{Conclusions}
\label{sec-7}

We have presented a new way to solve Poisson equation in the low-rank
tensor format using cross approximation in the frequency domain. The
idea is quite simple, but is very effective. Using the new solver we have
implemented the first low-rank solver for the two-dimensional Stokes
problem. The low-rank Stokes solver is based on the Uzawa-GMRES
approach, and is effective and robust in the experiments.  In the
numerical experiments only a two-dimensional case was considered,
since it is the ``worst case'' for the low-rank methods due to
typically high complexity with respect to the ranks. Even in this
setting, the low-rank solver outperforms the full-format
finite-difference solver starting from $n = 512$. 
We plan to implement a three-dimensional variant, and the complexity
reduction in three dimensions is expected to be much more visible. 
\bibliographystyle{siam}
\bibliography{final}

\end{document}